\documentclass{amsart}
\usepackage{amssymb, amsmath, amsfonts, amscd}

\input xypic

\theoremstyle{plain}
\newtheorem{theorem}{Theorem}[section]
\newtheorem{corollary}[theorem]{Corollary}
\newtheorem{lemma}[theorem]{Lemma}
\newtheorem{proposition}[theorem]{Proposition}

\newtheorem{example}[theorem]{Example}
\newtheorem{problem}[theorem]{Problem}

\theoremstyle{definition}
\newtheorem{definition}[theorem]{Definition}

\theoremstyle{remark}

\numberwithin{equation}{theorem}

\newcommand{\End}{\operatorname{End}}
\newcommand{\m}{\mathfrak{m}}
\newcommand{\F}{\mathcal{F}}
\newcommand{\G}{\mathbb{G}}
\newcommand{\I}{\mathcal{I}}

\renewcommand{\L}{\mathcal{L}}
\newcommand{\E}{\mathcal{E}}

\newcommand{\Q}{\mathcal{Q}} 
\renewcommand{\O}{\mathcal{O} }

\renewcommand{\P}{\mathbf{P} }
\renewcommand{\Pr}{\mathcal{J} }

\newcommand{\SL}{\operatorname{SL}}

\renewcommand{\H}{\operatorname{H} }
\newcommand{\Pic}{\operatorname{Pic} }

\renewcommand{\oe}{\overline{e}}
\newcommand{\R}{\operatorname{R} }
\newcommand{\U}{\operatorname{U}}

\renewcommand{\lg}{\mathfrak{g}}
\newcommand{\lp}{\mathfrak{p}} 

\renewcommand{\sl}{\mathfrak{sl}}
\newcommand{\n}{\mathfrak{n}_{-}}

\newcommand{\lh}{\mathfrak{h}}

\newcommand{\nplus}{\mathfrak{n}_{+}}
\newcommand{\nminus}{\mathfrak{n}_{-}}
\newcommand{\flag}{\mathbb{F}}
\newcommand{\dd}{\underline{d}}
\renewcommand{\ll}{\underline{l}}

\newcommand{\Z}{\mathbf{Z} }

\newcommand{\sym}{\operatorname{Sym} }

\newcommand{\one}{\mathbf{1}}
\newcommand{\spec}{\operatorname{Spec}}

\newcommand{\C}{F }

\begin{document}

\title{On jet bundles and generalized Verma modules II}

\author{Helge Maakestad}
\address{Institut Fourier, Grenoble}

\email{maakestd@fourier.ujf-grenoble.fr}
\keywords{linear algebraic group, jet bundle, vector bundle,
  flag variety, $P$-module, generalized Verma module, higher direct
  image, annihilator ideal}

\subjclass{14L30, 17B10, 14N15}

\date{Spring 2009}

\begin{abstract} Let $G$ be a semi simple linear algebraic group over
  a field of characteristic zero and let $V$ be a finite dimensional
  irreducible $G$-module with highest weight vector $v\in V$. Let
  $P\subseteq G$ be the parabolic subgroup fixing $v$. Let
  $\lg=Lie(G)$. We get a filtration $\U^\bullet(\lg)v: \U^k(\lg)v\subseteq V$ of
  $P$-modules for $1\leq k \leq N$. The aim of this paper is to use
  higher direct images of $G$-linearized sheaves, 
  filtrations of generalized Verma modules and annihilator ideals of highest weight
  vectors to give a natural basis for $\U^k(\lg)v$ and to compute its dimension.
  We also relate the filtration
  $\U^\bullet(\lg)v$ to $G$-linearized jet bundles on the flag variety $G/P$ for
  $G=\SL(E)$ where $E$ is a finite dimensional vector space.
\end{abstract}

\maketitle

\tableofcontents

\section{Introduction} 

Let $\C$ be a fixed field of characteristic zero and let $G$ be a semi
simple linear algebraic group over $\C$. Let $V$ be a finite
dimensional irreducible $G$-module with highest weight vector $v\in V$
and highest weight $\lambda$.  Let $L_v\subseteq V$ be the subspace
spanned by $v$.
Let $P\subseteq G$ be the parabolic
subgroup of elements fixing the subspace $L_v\subseteq V$. It follows
the quotient $G/P$ is a smooth projective variety of finite type over
$\C$. Let $\lg=Lie(G)$. We get a filtration of $V$ by $P$-modules

\begin{equation} \label{filtration1}
 \U^1(\lg)v\subseteq \U^2(\lg)v\subseteq \cdots
\subseteq \U^N(\lg)v=V 
\end{equation}
- the \emph{canonical filtration}.
Here $N=N(\lambda)$ is the minimal integer with $\U^N(\lg)v=V$. 

In a previous paper on this subject (see \cite{verma}) the filtration
\ref{filtration1} was studied in the case of $V=\H^0(\G(m,m+n),\O(d))^*$ on the
grassmannian $\SL(E)/P=\G(m,m+n)$. Here $P\subseteq \SL(E)$ is the
parabolic subgroup fixing an $m$-dimensional subspace in $E$.
There is an equivalence of categories between the category of
$\SL(E)$-linearized locally free sheaves on $\SL(E)/P$ and the
category of rational $P$-modules and the aim of this paper is to use
this equivalence to interpret the filtration \ref{filtration1}
in terms of $\SL(E)$-linearized locally free sheaves on $\SL(E)/P$.

We use higher direct images of $G$-linearized
sheaves, filtrations of generalized Verma
modules and annihilator ideals of highest weight vectors to 
answer the following questions for any parabolic subgroup  $P\subseteq
\SL(E)$ fixing a flag $E_\bullet$ in $E$ where $E$ is a finite dimensional vector space over $\C$:

\begin{align} 
&\label{main1}\text{Construct a basis for $\U^k(\lg)v$ as $\C$-vector space.}\\
&\label{main2}\text{Calculate the dimension  of $\U^k(\lg)v$.} \\
&\label{main3}\text{Interpret $\{\U^k(\lg)v\}_{k=0}^N$ in terms of geometric objects on $G/P$.}
\end{align}

The strategy of the proof is as follows. 

In section two of the paper we consider Question \ref{main1} and \ref{main2}.
Let \[ E_\bullet: 0 \neq E_1 \subseteq E_2\subseteq \cdots \subseteq E_k
\subseteq
E_{k+1}=E \]
be a flag in the vector space $E$ of type
$\underline{d}=(d_1,d_2,..,d_k,d_{k+1})$. This means
$dim_\C(E_i)=d_1+\cdots
+d_i$. Let $n_i=d_1+\cdots +d_i$ for $i=1,..,k+1$. Let $P\subseteq
G=\SL(E)$ be the parabolic subgroup fixing the flag $E_\bullet$.
The quotient $G/P=\flag(\dd,E)$ is the \emph{flag variety of type
  $\dd$} parametrizing flags in $E$. Let $V$ have highest weight
\[ \lambda=\sum_{i=1}^k l_i(L_1+\cdots + L_{n_i}).\]
Let $\L(\ll)\in \Pic^G(G/P)=\mathbb{Z}^k$ be the line bundle
corresponding to $\ll=(l_1,..,l_k)$. The $G$-module
$\H^0(G/P,\L(\ll))^*$ has highest weight $\lambda$ hence there is an isomorphism of
$G$-modules $V\cong \H^0(G/P,\L(\ll))^*$ giving a
geometric construction of the $G$-module $V$. 
The highest weight vector $v$ of $V$ has a geometric construction. It
is the vector defined by $v:\H^0(G/P,\L(\ll))\rightarrow \C$ and $v(s)=s(\overline{e})$.
There is an inclusion of $G$-modules
\begin{equation}\label{inclusion}
\H^0(G/P,\L(\ll))^* \subseteq \sym^{l_1}(\wedge^{d_1}E)\otimes
\cdots \otimes \sym^{l_k}(\wedge^{d_k}E)
\end{equation}
and the highest weight
vector $v\in V=\H^0(G/P,\L(\ll))^*$ is via the inclusion
\ref{inclusion} described explicitly as follows:
\begin{equation} \label{highest}
 v=\sym^{l_1}(\wedge^{d_1}E_1)\otimes \cdots \otimes
\sym^{l_k}(\wedge^{d_k}E_k)\in \H^0(G/P,\L(\ll))^*.
\end{equation}
We get an exact sequence of $G$-modules
\[ 0\rightarrow ann(v,\lambda)\otimes L_v\rightarrow \U(g)\otimes L_v
\rightarrow \H^0(G/P,\L(\ll))^* \rightarrow 0\]
and an exact sequence of $P$-modules
\begin{equation}\label{ex1}
 0\rightarrow ann^k(v,\lambda)\otimes L_v \rightarrow
\U^k(\lg)\otimes L_v \rightarrow \U^k(\lg)v\rightarrow 0. 
\end{equation}
Here $ann(v,\lambda)\subseteq \U(\lg)$ is the left \emph{annihilator ideal}
of $v\in V$. The left $G$-module $\U(\lg)\otimes L_v$ is a \emph{generalized Verma
  module} and $\U^\bullet(\lg)\otimes L_v\subseteq \U(\lg)\otimes L_v$
is the canonical filtration of $\U(\lg)\otimes L_v$. 
Let $\lp=Lie(P)$ be the Lie algebra fixing $v\in V$.
We use the explicit description of $v\in V=\H^0(G/P,\L(\ll))^*$ given
in \ref{highest}, the exact sequence
\ref{ex1} and properties of the universal enveloping algebra $\U(\lg)$ 
to prove the following Theorem:

\begin{theorem} Let $1\leq k \leq min\{l_i+1\}$. There is an equality of
vector spaces
\[ W^k(v,\lambda)=\U^k(\n) .\]
Here $W^k(v,\lambda)$ is a complement of $ann^k(v,\lambda)$ in
$\U^k(\lg)$ and $\n \subseteq \lg$ is a sub Lie algebra with $\lg=\n\oplus \lp$.
\end{theorem}
\begin{proof} See Theorem \ref{maintheorem1}.
\end{proof}

From Theorem \ref{maintheorem1} we can prove Corollary \ref{corrmain1}
answering Question \ref{main1} and 
Corollary \ref{corrmain2} answering Question \ref{main2}.

In the final section we consider Question \ref{main3}. 
The vector space $\U^k(\lg)v$ is a
$P$-module hence it corresponds to an $\SL(E)$-linearized locally free
sheaf $\Pr^k$ on $\SL(E)/P$. We make the equivalence explicit and give a geometric
construction of $\Pr^k$ in terms of $\SL(E)$-linearized jet bundles on 
$\SL(E)/P$. We use a
vanishing theorem from an earlier paper on the subject (see
\cite{verma}), Kunneth formulas and general properties of jet bundles
on products to construct an exact sequence of $P$-modules
\[ 0\rightarrow \Pr^k_{\flag}(\L(\ll)(\overline{e})^*\rightarrow
\H^0(\flag,\L(\ll))^* \rightarrow^\phi
\H^0(\flag,\m^{k+1}\L(\ll))^*\rightarrow 0 \]
when $1\leq k \leq min\{l_i+1\}$. We get an injection of $P$-modules
\[ \Pr^k_{\flag}(\L(\ll))(\overline{e})^*\subseteq
V=\H^0(\flag(\dd,E),\L(\ll))^*.\]
Here we use the fact we can give a geometric construction of the
$\SL(E)$-module $V$. Then we prove the main result of the paper:

\begin{theorem} Let $1\leq k \leq min\{l_i+1\}$. There is an
  isomorphism of $P$-modules
\[ \U^k(\lg)v\cong \Pr^k_{\flag}(\L(\ll)(\overline{e})^* .\]
\end{theorem}
\begin{proof} See Theorem \ref{maintheorem2}.
\end{proof}

Theorem \ref{maintheorem2} answer Question \ref{main3} giving a
geometric interpretation of the filtration $\U^\bullet(\lg)v \subseteq
V$ in terms of $\SL(E)$-linearized jet bundles on $\SL(E)/P$.

The motivation for the study of the jet bundle $\Pr^k(\L(\ll))$ is
partly its relationship with the discriminant $D^k(\L(\ll))$ of the
line bundle $\L(\ll)$. Assume $\ll=(l_1,..,l_k)\in \Z^k$ with $l_i\geq
1$ for all $i$. It follows by the results of this paper the $k$'th
Taylor map
\[ T^k:\H^0(G/P, \L(\ll))\otimes \O_{G/P}\rightarrow \Pr^k(\L(\ll)) \]
is surjective when $1\leq k \leq min\{l_i+1\}$. We get an exact
sequence of locally free sheaves
\[ 0\rightarrow \Q_{\L(\ll)}\rightarrow \H^0(G/P, \L(\ll))\otimes
\O_{G/P}\rightarrow \Pr^k(\L(\ll)) \rightarrow 0\]
of $\O_{G/P}$-modules. Dualize this sequence to get an exact sequence 
\[ 0\rightarrow \Pr^k(\L(\ll))^* \rightarrow \H^0(G/P,
\L(\ll))^*\otimes \O_{G/P} \rightarrow \Q_{\L(\ll)}^* \rightarrow 0.\]
Take relative projective space bundle to get a closed immersion of
schemes
\[ \P(\Q^*_{\L(\ll)})\subseteq \P(W^*)\times G/P \]
where
\[ W=\H^0(G/P, \L(\ll)).\]
There is a projection map
\[p:\P(W^*)\times G/P \rightarrow \P(W^*) \]
and by the results of \cite{maa12} it follows the direct image scheme
$p(\P(\Q^*_{\L(\ll)}))$ equals the discriminant $D^k(\L(\ll))$ of the
line bundle $\L(\ll)$. There is on $\P(W^*)$ the tautological sequence
\[ 0\rightarrow \O(-1)\rightarrow W\otimes \O_{\P(W^*)} .\]
pull this and the Taylor map back to $Y=\P(W^*)\times G/P$ to get the
composed map
\[ \phi:\O(-1)_Y\rightarrow W\otimes
\O_{G/P}\rightarrow^{T^k} \Pr^k(\L(\ll))_Y.\]
It follows by the results of \cite{maa12}  the scheme theoretic image
of the zero scheme
$p(Z(\phi))$ equals the discriminant $D^k(\L(\ll))$ as subscheme of
$Y$. When the ideal sheaf of $Z(\phi)$ is locally generated by a
regular sequence we get a Koszul complex of locally free sheaves
\[ 0\rightarrow \O(-r)_Y\otimes \wedge^r \Pr^k(\L(\ll))^*_Y\rightarrow
\cdots \rightarrow \O(-1)_Y\otimes \Pr^k(\L(\ll))^*_Y\rightarrow \O_Y
\rightarrow \O_{Z(\phi)}\rightarrow 0 \]
which is a resolution of the ideal sheaf of $Z(\phi)$. When we push
this complex down to $\P(W^*)$ we get a double complex with terms as
follows:
\[ \R^ip_*(\O(-j)_Y\otimes \wedge^j\Pr^k(\L(\ll))^*_Y)=\O(-j)\otimes
\H^i(G/P, \wedge^j\Pr^k(\L(\ll)^*) .\]
By \cite{maa12}, Theorem 5.2 we know knowledge on the $P$-module structure of the fiber
$\Pr^k(\L(\ll))(e)^*$ will give information on
$\H^i(G/P,\wedge^j\Pr^k(\L(\ll))^*)$. Hence we may determine if the
double complex
\[ \O(-j)\otimes \H^i(G/P, \wedge^j\Pr^k(\L(\ll)^*) .\]
can be used to construct a resolution of the ideal sheaf $\I$ of the
discriminant $D^k(\L(\ll))$.

\section{Character ideals and annihilator ideals}

Let in this section $E$ be a fixed $N$-dimensional vector space over
$\C$.
Let
\[ E_\bullet: 0 \neq E_1 \subseteq E_2\subseteq \cdots \subseteq E_k \subseteq
E_{k+1}=E \]
be a flag in the vector space $E$ of type
$\underline{d}=(d_1,d_2,..,d_k,d_{k+1})$. This means $dim_\C(E_i)=d_1+\cdots
+d_i$. Let $n_i=d_1+\cdots +d_i$ for $i=1,..,k+1$. Let $P\subseteq
\SL(E)$ be the parabolic subgroup fixing the flag $E_\bullet$. This
means $g(E_i)\subseteq E_i$ for all $g\in P$ and $i=1,..,k$.
The quotient $\flag(\dd,E)=\SL(E)/P$ is a smooth projective variety of finite type
over $\C$ - the \emph{flag variety of type $\dd$}. It is a
\emph{geometric quotient} in the sense of \cite{mumford}
and it is the parameter space parametrizing flags of type $\dd$ in
$E$. This means each point $x\in \flag(\dd,E)$ with coefficients in $\C$ corresponds to a
unique flag
\[ E_1^x\subseteq E_2^x \subseteq \cdots \subseteq E_k^x \subseteq E
\]
of type $\dd$.

Let
$\G_i=\G(n_i,N)$ be the grassmannian of $n_i$-planes in an $N$-dimensional $\C$-vector space.
There is a closed immersion - the \emph{generalized Plucker embedding}
\[ i: \flag(\dd, E) \rightarrow \G_1\times \cdots \times \G_k
\]
defined by
\[ i([E_1\subseteq E_2 \subseteq \cdots \subseteq E_k ])=([E_1],[E_2],..,[E_k]).\]
Let $\ll=(l_1,..,l_k)\in \mathbb{Z}^k$ be a $k$-tuple of integers. We
get a line bundle 
\[  \L(\ll)=i^*\O(l_1)\otimes \cdots \otimes \O(l_k) \]
and the line bundle $\L(\ll)$ has by \cite{mumford} a unique $\SL(E)$-linearization. We
get an isomorphism
\[ \Pic^{\SL(E)}(\flag(\dd,E))\cong \mathbb{Z}^k \]
of groups. If $l_i\geq 1$ for $i=1,..,k$ it follows $\L(\ll)$ is very
ample. It follows there is a closed immersion
\[ j:\flag(\dd,E)\rightarrow \mathbb{P}^M \]
with $M>>0$ and  $j^*\O(1)=\L(\ll)$.

Let $P\subseteq G$ be the subgroup consisting
  of matrices $g$ with determinant one of the following type:
\[
g=\begin{pmatrix} A_1 & *      &  \cdots    & *      &     * \\
                 0  &   A_2  &  \cdots    & *      &     * \\
             \vdots &  \vdots&  \cdots    & \vdots & \vdots \\
                 0  &   0    &  \cdots    & A_k    &  *     \\
                 0  &   0    &  \cdots    & 0      & A_{k+1}
\end{pmatrix}
\]
where $A_i$ is an $d_i \times d_i$-matrix with coefficients in $\C$.
It follows $P$ is the group of elements $g\in \SL(E)$ fixing the flag $E_\bullet$.
Note: We can define the subgroup $P$ using points with values in $\C$
since $\C$ has characteristic zero and all group schemes in
characteristic zero are smooth.
It follows $P$ is a parabolic group and it follows the quotient
$G/P$ equals $\flag(\dd,E)$. Let $\ll=(l_1,..,l_k)\in \mathbb{Z}^k$ and let
$\L(\ll)\in \Pic^G(G/P)$ be the line bundle defined above. 
There is a unique $P$-stable vector $v$ in $\H^0(G/P,\L(\ll))^*$ defined as follows:
\[ v:\H^0(G/P,\L(\ll))\rightarrow \C \]
with
\[ v(s)=s(\overline{e}) \]
where $\overline{e}\in G/P$ is the class of the identity element.

\begin{lemma} \label{highestweight}
The vector $v$ is a highest weight vector for $\H^0(G/P,\L(\ll))^*$
with highest weight
\[ \lambda =\sum_{i=1}^kl_i \omega_{n_i} \]
\end{lemma}
\begin{proof}
This is left to the reader as an exercise.
\end{proof}

Let $V$ be an arbitrary finite dimensional irreducible $G$-module with
highest
weight vector $v\in V$ and highest weight $\lambda$.
Let $\omega_j=L_1+\cdots +L_j$ for $j=1,..,N-1$ be the fundamental
weights for $G=\SL(E)$. It follows $\lambda=\sum_{i=1}^k
l_i\omega_{n_i}$ with $l_i\geq 0$ for all $i$. There is the following well known result:

\begin{theorem} \label{BWB}There is a parabolic subgroup $P\subseteq G$ and
  a linebundle $\L(\ll)\in \Pic^G(G/P)$ with an isomorphism $V\cong
  \H^0(G/P,\L(\ll))^*$ as $G$-modules.
\end{theorem}
\begin{proof} Let $P\subseteq G$ be the above defined parabolic
  subgroup and let $\L(\ll)$ be the line bundle with $\ll=(l_1,..,l_k)\in\mathbf{Z}^k$
It follows by the \emph{Borel-Weil-Bott Theorem} (see \cite{akhiezer})
that $\H^0(G/P,\L(\ll))^*$ is an irreducible $G$-module.
By Lemma \ref{highestweight} it follows the $P$-stable vector $v\in
\H^0(G/P,\L(\ll))^*$
has highest weight $\lambda$ with
\[ \lambda= \sum_{i=1}^k l_i(tr(A_1)+\cdots +tr(A_i))=\sum_{i=1}^kl_i
\omega_{n_i} \]
hence there is an isomorphism
\[ V\cong \H^0(G/P,\L(\ll))^* \]
of $\SL(E)$-modules and the Theorem is proved.
\end{proof}

Hence Theorem \ref{BWB} gives a geometric construction of all finite
dimensional irreducible representations of $G=\SL(E)$: All
irreducible finite dimensional $G$-modules may be realized as duals of
global
sections of $G$-linearized line bundles on $G/P$ for some parabolic
subgroup $P \subseteq G$. Note: The subgroup $P$ is not unique.

Let $L_v$ be the subspace spanned by $v$.
The group $P\subseteq \SL(E)$ is the subgroup of elements $g\in
\SL(E)$ stabilizing the space $L_v$ defined by the highest weight vector $v\in
V=\H^0(G/P,\L(\ll))^*$. Let $\lg=Lie(\SL(E))$ and let
$\U^k(\lg)\subseteq \U(\lg)$ be the canonical filtration of the
universal enveloping algebra of $\lg$. The vector space $\U^k(\lg) $
is a $G$-module via the adjoint representation. It follows $\U^k(\lg)$
is a $P$-module. There is a surjective map of $G$-modules
\[ \phi:\U(\lg)\otimes L_v\rightarrow \H^0(\G/P,\L(\ll))^* \]
defined by
\[\phi( x\otimes v) = x(v).\] The $G$-module $\U(\lg)\otimes L_v$
is a \emph{generalized Verma module} and the map $\phi$ realize the
$G$-module $\H^0(G/P,\L(\ll))^*$ as a quotient of $\U(\lg)\otimes
L_v$. It is a fact that any finite dimensional irreducible $G$-module
may be realized as a quotient of a generalized Verma module.
The map $\phi$ induces a surjection of $P$-modules
\[ \U^k(\lg)\otimes L_v\rightarrow \U^k(\lg)v  \]
where $\U^k(\lg)v$ is the vector space of elements $x(v)$ with $x\in
\U^k(\lg)$.
We get a filtration of $V$ by $P$-modules:
\[ \U^\bullet(\lg)v: \U^1(\lg)v \subseteq \cdots \U^k(\lg)v \subseteq V.\]

\begin{definition} Let the filtration $\U^\bullet(\lg)v\subseteq V$ 
be the \emph{canonical filtration} of $V$.
\end{definition}

Since the $P$-module $\U^k(\lg)v$ only depends on the vector space
$L_v$ defined by the highest weight vector $v\in V$ we have defined
for an arbitrary irreducible $\SL(E)$-module $V$ a canonical
filtration $\U^\bullet(\lg)v\subseteq V$. 

Note: This notion is well defined for an arbitrary irreducible finite
dimensional representation of an arbitrary semi simple Lie algebra.

\begin{definition}
Let $ann(v,\lambda )\subseteq \U(\lg)$ be the left
\emph{annihilator ideal} of $v\in V$ and let $ann^k(v,\lambda)\subseteq
ann(v,\lambda)$ be its canonical filtration. 
\end{definition}

We get an exact sequence of $G$-modules 
\begin{equation} \label{vermaquotient}
 0\rightarrow ann(v,\lambda)\otimes L_v\rightarrow \U(\lg)\otimes
L_v \rightarrow \H^0(G/P,\L(\ll))^* \rightarrow 0
\end{equation}

and an exact sequence of $P$-modules
\begin{equation} 
\label{exactsequence}  0\rightarrow ann^k(v,\lambda)\otimes L_v\rightarrow
\U^k(\lg)\otimes L_v \rightarrow \U^k(\lg)v \rightarrow 0 
\end{equation}

where the rightmost map is the obvious action map. 
The sequence \ref{exactsequence} describe the terms $\U^k(\lg)v$ in the canonical filtration
\[\U^\bullet(\lg)v: \{v\}\subseteq \U^1(\lg)v \subseteq
\U^2(\lg)v\subseteq \cdots \subseteq \U^N(\lg)v=V \]
of $V=\H^0(G/P,\L(\ll))^*$ as quotients of the terms in the filtration $\U^\bullet(\lg)\otimes L_v$
of the generalized Verma module $\U(\lg)\otimes L_v$. Here $N=N(\lambda)$ is the minimal integer
with the property that $\U^N(\lg)v=V$.
The aim of this section is to answer the questions posed in
\ref{main1} and \ref{main2} using the exact sequence \ref{exactsequence} and
properties of the universal enveloping algebra.

Let $\lp=Lie(P)$ be the Lie algebra of $P$. It follows $\lp\subseteq
\lg=\sl(E)$ is the Lie algebra of matrices $x\in \sl(E)$ of the form
\[
x=\begin{pmatrix} A_1 & *      &  \cdots    & *      &     * \\
                 0  &   A_2  &  \cdots    & *      &     * \\
             \vdots &  \vdots&  \cdots    & \vdots & \vdots \\
                 0  &   0    &  \cdots    & A_k    &  *     \\
                 0  &   0    &  \cdots    & 0      & A_{k+1}
\end{pmatrix}.
\]
Here $A_i$ is an $d_i\times d_i$-matrix with coefficients in $\C$ and $tr(x)=0$. 
We aim to calculate a vector space
$W^k(v,\lambda)\subseteq \U^k(\lg)$ with the property that
\[ \U^k(\lg)=W^k(v,\lambda)\oplus ann^k(v,\lambda)  .\]
The line $(v)$defined by $v$ is $\lp$-stable and we get a character
\[ \rho:\lp\rightarrow \End_\C(v) \]
defined by
\[ \rho(x)(v)=xv.\]
Since $\End_\C(v)\cong \C$ we get a map
\[ \rho:\lp \rightarrow \C .\]
One checks that
\[ \rho(x)=\sum_{i=1}^kl_i(tr(A_1)+\cdots +tr(A_i)) .\]
\begin{definition} Let $char(\rho)=\U(\lg)\{x-\rho(x): x\in \lp\}$ be
  the left \emph{character ideal} of $\rho$. Let 
$char^k(\rho)=char(\rho)\cap \U^k(\lg)$ be its canonical filtration.
\end{definition}
Let $\n\subseteq \lg$ be the complement of $\lp\subseteq \lg$. It is
a  sub Lie algebra.

Let $x\in \lp$ be the following matrix:
\[
x=\begin{pmatrix} 0 & 0      &  \cdots    & 0      &     0 \\
                 0  &   0  &  \cdots    & 0      &     0 \\
             \vdots &  \vdots&  \cdots    & \vdots & \vdots \\
                 0  &   0    &   \vdots   & 0      & 0      \\
                 0  &   0    &  \cdots    & A_k    &  0     \\
                 0  &   0    &  \cdots    & 0      & A_{k+1}
\end{pmatrix}
\]
where $A_k$ is the $d_k\times d_k$  matrix
\[
A_k=\begin{pmatrix} 0 & 0      &  \cdots    & 0      &     0 \\
                 0  &   0  &  \cdots    & 0      &     0 \\
             \vdots &  \vdots&  \cdots    & \vdots & \vdots \\
                 0  &   0    &  \cdots    & 0    &  0     \\
                 0  &   0    &  \cdots    & 0      & 1
\end{pmatrix}
\]
and $A_{k+1}$ is the $d_{k+1}\times d_{k+1}$-matrix
\[
A_{k+1}=\begin{pmatrix} 0 & 0      &  \cdots    & 0      &     0 \\
                 0  &   0  &  \cdots    & 0      &     0 \\
             \vdots &  \vdots&  \cdots    & \vdots & \vdots \\
                 0  &   0    &  \cdots    & 0    &  0     \\
                 0  &   0    &  \cdots    & 0      & -1
\end{pmatrix}.
\]

Let $x_{n_j}$ with $n_{i-1}+1\leq j \leq n_i$ with $i=1,..,k$ be the following matrix:
\[
x_{n_i}=\begin{pmatrix} 0 & 0      &  \cdots    & 0      &     0 \\
                 0  &   0  &  \cdots    & 0      &     0 \\
             \vdots &  \vdots&  A_i    & \vdots & \vdots \\
                 0  &   0    &  \cdots    & 0    &  0     \\
                 0  &   0    &  \cdots    & 0      & A_{k+1}
\end{pmatrix}
\]
Where $A_i$ is the matrix with zeros everywhere and $1$ on the $j$'th
place on the diagonal. The matrix $A_{k+1}$ has zeros everywhere and
$-1$ in the lower right corner. 

Let $\lp\subseteq \lg$ be the stabilizer Lie algebra of $v$ and let $\lp_v\subseteq \lp$ be the
isotropy Lie algebra of $v\in \H^0(G/P,\L(\ll))^*$. Let $\lp =\lp_0
\oplus \lp_D$ where
$\lp_0$ is the subspace of matrices with zeros on the diagonal and
$\lp_D\subseteq \lp$ is the subspace of diagonal matrices. 
It follows $\lp_0\subseteq \lp_v$.
\begin{lemma} The set  
\[ \{x_{n_j}: i=1,..,k; n_{i-1}+1\leq j \leq n_i\} \]
is a basis for the vector space $\lp_D$. 
\end{lemma}
\begin{proof} The proof is left to the reader as an exercise.
\end{proof}

\begin{proposition} \label{element}Let $v\in \H^0(G/P,\L(\ll))^*$ be
  the unique highest weight vector. It follows $\lp=\lp_v\oplus (x)$ where $x$ is the matrix
  defined above. Furthermore $x(v)=l_kv$.
\end{proposition}
\begin{proof} 
Let for any matrix $z\in \lp$
\[
z=\begin{pmatrix} A_1 & *      &  \cdots    & *      &     * \\
                 0  &   A_2  &  \cdots    & *      &     * \\
             \vdots &  \vdots&  \cdots    & \vdots & \vdots \\
                 0  &   0    &  \cdots    & A_k    &  *     \\
                 0  &   0    &  \cdots    & 0      & A_{k+1}
\end{pmatrix}.
\]
$B_i$ be the associated submatrix
\[
B_i=\begin{pmatrix} A_1 & *      &  \cdots    & *      &     * \\
                 0  &   A_2  &  \cdots    & *      &     * \\
             \vdots &  \vdots&  \cdots    & \vdots & \vdots \\
                 0  &   0    &  \cdots    & A_{i-1}    &  *     \\
                 0  &   0    &  \cdots    & 0      & A_{i}
\end{pmatrix}.
\]

A matrix $z\in \lp$ is in $\lp_v$ if and only if 
\[ \rho(z)=\sum_{i=1}^k l_itr(B_i)=0.\]
Hence $z\in \lp-\lp_v$if and only if $\rho(z)\neq 0$. 
Write $\lp=\lp_0\oplus \lp_D$. In the discussion preceeding the
Proposition we constructed a basis $x_{n_j}$ for $\lp_D$ with $n_{i-1}+1 \leq j
\leq n_i$ for $i=1,..,k$. By definition it follows $x_j\in \lp_v$ for
$n_k+1\leq j \leq n_{k+1}$. Hence an element
\[ y=\sum_{s=1}^{n_k}a_sx_s \]
is in $\lp_v$ if and only if
\[ \sum_{t=1}^k l_ttr(B_t)=0.\]
This is if and only if there is an equation
\[ a_{n_k}=f(a_1,..,a_{n_k-1}).\]
Hence $y\notin \lp_v$ if and only if 
\[ a_{n_k}\neq f(a_1,..,a_{n_k-1}) \]
and we check that the only element in the above constructed basis for
$\lp_D$ satisfying this condition is the element $x$
defined above, and the first claim of the Proposition follows.
One checks the second claim of the Proposition by calculation and the
Proposition is proved.
\end{proof}

\begin{lemma} \label{sum} Let $v\in \sym^{k+1}(\lg)\subseteq \U^{k+1}(\lg)$ be an
  element. We may write $v=v_1+v_2$ with $v_1\in \sym^{k+1}(\n)$ and 
$v_2\in  \U^k(\lg)\{y-\rho(y): y\in \lp\}$.
\end{lemma}
\begin{proof} The proof is left to the reader as an exercise.
\end{proof}

The element $x$ from Proposition \ref{element} depends on the
decomposition $\lp=\lp_v\oplus (x)$ but this fact will not be
important in what follows.

\begin{proposition} \label{decomposition} There is for all $k\geq 1$ an equality
\[ \U^k(\lg)=\U^k(\n)\oplus char^k(\rho) \]
of vector spaces.
\end{proposition}
\begin{proof}  One checks that there is an equality of vector spaces
\[ char^1(\rho)=\{x-l_k\one, y : y\in \lp_v\}.\] 
Using the Poincare-Birkhoff-Witt Theorem one checks there is an
equality of vector spaces
\[ char^k(\rho)=\U^{k-1}(\lg)\{y-\rho(y):y\in \lp\}.\]
We prove the claim in the Proposition using induction on $k$. We first
check it for $k=1$. We get
\[ \U^1(\lg)=\one\oplus \lg=\one \oplus \n\oplus (x)\oplus \lp_v=\]
\[ \one \oplus \n \oplus (x-l_k\one)\oplus \lp_v = \U^1(\n)\oplus
\{x-l_k\one,y:y\in \lp\}=\]
\[\U^1(\n)\oplus \{y-\rho(y):y\in \lp\} =\U^1(\n)\oplus
char^1(\rho),\]
and the claim of the Proposition is proved for $k=1$.
Assume the claim is true for $k$:
\[ \U^k(\lg)=\U^k(\n)\oplus char^k(\rho).\]
Using the symmetrization map we may identify 
\[ \U^k(\lg)=\oplus_{i=0}^k\sym(\lg) \]
where $\sym(\lg)$ is the $i$'th symmetric power of $\lg$ with the
adjoint representation. It follows there is a isomorphism
\[ \U^{k+1}(\lg)\cong \U^k(\lg)\oplus \sym^{k+1}(\lg) \]
of $\lg$-modules. Since the symmetrization map is an isomorphism of
vector spaces we may identify $\sym^{k+1}(\lg)$ with its image in
$\U(\lg)$. All calculations in what follows are done inside $\U(\lg)$
via the symmetrization map.
Clearly there is an inclusion
\[ \U^{k+1}(\n)\oplus \U^k(\lg)\{y-\rho(y):y\in \lp\}\subseteq
\U^{k+1}(\lg).\]
We prove the reverse inclusion. Write
\[ \U^{k+1}(\lg)=\U^k(\lg)\oplus\sym^{k+1}(\lg)=\]
\[\U^k(\n)\oplus\U^{k-1}(\lg)\{y-\rho(y):y\in \lp\}\oplus
\sym^{k+1}(\lg).\]
Let $v\in \sym^{k+1}(\lg)$. From Lemma \ref{sum}  one may write
$v=v_1+v_2$ with $v_1\in \sym^{k+1}(\n)$ and $v_2\in
\U^k(\lg)\{y-\rho(y):y\in \lp\}$. It follows
\[\U^{k+1}(\lg)=\U^{k+1}(\n)\oplus \U^k(\lg)\{y-\rho(y):y\in
\lp\}=\U^{k+1}(\n)\oplus char^{k+1}(\rho) \]
and the claim of the Proposition is proved.
\end{proof}

Consider the Plucker embedding
\[ i:\flag(\dd,E)\rightarrow \G_1\times \cdots \times \G_k \subseteq \mathbb{P}^M\]
defined in the beginning of this section. The flag variety is
projectively normal hence there is an injection of vector spaces
\[ \H^0(\flag(\dd,E),\L(\ll))^* \subseteq
\sym^{l_1}(\wedge^{n_1}E)\otimes \cdots \otimes
\sym^{l_k}(\wedge^{n_k}E) .\]
There is a $P$-stable line 
\[ \sym^{l_1}(\wedge^{n_1}E_1)\otimes \cdots \otimes
\sym^{l_k}(\wedge^{n_k}E_k) \subseteq \sym^{l_1}(\wedge^{n_1}E)\otimes
\cdots \otimes \sym^{l_k}(\wedge^{n_k}E) .\]

\begin{lemma} Let $v\in \H^0(\flag(\dd,E),\L(\ll))^*$ be the highest
  weight vector. There is an equality
\[ v= \sym^{l_1}(\wedge^{n_1}E_1)\otimes \cdots \otimes \sym^{l_k}(\wedge^{n_k}E_k) .\]
\end{lemma}
\begin{proof}
This is left to the reader as an exercise.
\end{proof}

Let $v_i=\wedge^{n_i}E_i$ for $i=1,..,k$. We write $v=v_1^{l_1}\otimes
\cdots \otimes v_k^{l_k}$.

We use the notation of
\cite{dixmier} Chapter $7$. Let $P$ be the dominant weights of $\lg=\sl(E)$ and
let $B$ be a basis for the roots of $\lg$. Let $\rho$ be the character
associated to $v$. It follows for $x\in \lp$ we have
\[ \rho(x)=\sum_{i=1}^k l_i(tr(A_1)+\cdots +tr(A_i)) .\]
Let $\lg=\lg_{-}\oplus \lh \oplus \lg_{+}$ be a Cartan decomposition
of $\lg$. Let $\nplus =\lg_{+}\oplus \lh$.
Let $I(v)$ be the left ideal in $\U(\lg)$ defined by
\[ I(v)=\U(\lg)\nplus +\sum_{x\in \lh}\U(\lg)(x-\lambda(x)) .\]
By \cite{dixmier} Proposition 7.2.7 it follows
\[ ann(v,\lambda)=I(v)+\sum_{\beta\in
  B}\U(\nminus)X^{m_\beta}_{-\beta}.\]
Let $I^k(v)=I(v)\cap \U^k(\lg)$ be the canonical filtration of $I(v)$.
Let $B=\{L_i-L_{i+1}\}_{i=1,..,N}$. Let $\beta_i=L_i-L_{i+1}$ and let
$\lg^{\beta_i}=\C (E_{i,i+1})$.
It follows $\lg^{-\beta_i}=\C (E_{i+1,i})$. We have by definition
$X_{-\beta_i}=E_{i,i+1}$ and since $[E_{i,j},E_{j,i}]=E_{i,i}-E_{j,j}$
and $0\neq H_{\beta_i}\in [\lg^{\beta_i},\lg^{-\beta_i}]$
it follows
\[ H_{\beta_i}=E_{i,i}-E_{i+1,i+1}.\]
By definition we have
\[ m_{\beta_i}=\lambda(H_{\beta_i})+1.\]

\begin{lemma}
\begin{align}
m_{\beta_i}=&\label{mbeta1} l_j+1 \text{ if $i=n_j$} \\
m_{\beta_i}=&\label{mbeta2} 1\text{ if $i\neq n_j$}
\end{align} 
\end{lemma}
\begin{proof} The proof is left to the reader as an exercise.
\end{proof}

Let $K^k(v)$  be the following vector space:
\[ K^k(v)=(\sum_{\beta\in B}\U(\nminus)X^{m_{\beta}}_{-\beta})\cap
\U^k(\lg).\]
It follows
\[ K^k(v)=\sum_{i\neq
  n_j}\U^{k-1}(\nminus)X_{-\beta_i}+\sum_{i=n_j}\U^{k-l_j-1}(\nminus)X^{l_j+1}_{-\beta_{n_j}}.\]
If $1\leq k \leq min\{l_i+1\}$ it follows
\[ K^k(v)=\sum_{i\neq n_j}\U^{k-1}(\nminus)X_{-\beta_i}.\]
We have
\[ann^k(v,\lambda)=I^k(v)+K^k(v) .\]

\begin{theorem} \label{maintheorem1}Let $1\leq k \leq min\{l_i+1\}$ be an integer. The following holds:
\[ W^k(v,\lambda)=\U^k(\n) .\]
\end{theorem}
\begin{proof} By definition it follows $char^k(\rho)\subseteq
  ann^k(v,\lambda)$ for all $k\geq 1$.
We want
to prove the reverse inclusion 
\[ ann^k(v,\lambda)\subseteq char^k(\rho) \]
in the case when $1\leq k \leq min\{l_i+1\}$.

There is an inclusion $I^k(v)\subseteq char^k(\rho)$ for all $k\geq
1$. When $1\leq k\leq min\{l_i+1\}$ we get $K^k(v)\subseteq
char^k(\rho)$ and it follows
\[ ann^k(v,\lambda)\subseteq char^k(\rho).\]
From this we deduce an equality 
\[ ann^k(v,\lambda)=char^k(\rho) \]
when $1\leq k\leq min\{l_i+1\}$.
By Proposition \ref{decomposition} the following holds:
\[\U^k(\lg)=\U^k(\n)\oplus char^k(\rho) \]
It follows $W^k(v,\lambda)=\U^k(\n)$  and the claim of the Theorem is proved.
\end{proof}

\begin{corollary} \label{equalityv} Let $1\leq k \leq min\{l_i+1\}$. There is an
  equality of vector spaces
\[ \U^k(\lg)=\U^k(\n)\oplus ann^k(v,\lambda).\]
\end{corollary}
\begin{proof} This follows from Theorem \ref{maintheorem1}.

\end{proof}

\begin{corollary} \label{corrmain1} Let $v_1,..,v_D$ be a basis for
  $\n\subseteq \lg$ and let $1\leq k \leq min\{l_i+1\}$.
It follows the set
\[ \{v_1^{a_1}\cdots v_D^{a_D}(v): 0\leq \sum_i a_i \leq k\} \]
is a basis for $\U^k(\lg)v.$ 
\end{corollary}
\begin{proof} There is  by Corollary \ref{equalityv} an equality
\[ \U^k(\lg)=\U^k(\n)\oplus ann^k(v,\lambda) .\]
It follows from this there is an isomorphism of vector spaces
\[ \U^k(\n)\otimes L_v\rightarrow \U^k(\lg)v .\]
From this isomorphism and the Poincare-Birkhoff-Witt Theorem 
the claim of the Corollary follows.
\end{proof}

Let $D=\sum_{1\leq i < j\leq k+1}d_id_j$ It follows $dim_\C(\n)=D$. 

\begin{corollary} \label{corrmain2} Let $1\leq k \leq min\{l_i+1\}$. The following holds:
\[dim_\C(\U^k(\lg)v= \binom{D+k}{D}.\]
\end{corollary}
\begin{proof} There is by Theorem \ref{maintheorem1} an isomorphism of
  vector spaces $\U^k(\n)\cong \U^k(\lg)v$. It follows
\[dim_\C(\U^k(\lg)v)=dim_\C(\U^k(\n))=dim_\C(\sym^k(\n\oplus
\one))=\binom{D+k}{D}\]
and the claim of the Corollary is proved.
\end{proof}

Corollary \ref{corrmain1} and \ref{corrmain2} answer the questions
\ref{main1} and \ref{main2} posed in the introduction of the paper.

\section{Filtrations of $\SL(E)$-modules and jet bundles}

In this section we relate filtration $\U^k(\lg)v\subseteq V$ studied in the
previous section to the jet bundle $\Pr^k_{\flag}(\L(\ll))$ of the line bundle
$\L(\ll)\in \Pic^G(G/P)$ with $\H^0(G/P,\L(\ll))^*=V$.
Recall: The vector space $\U^k(\lg)v$ is a
$P$-module hence it corresponds to an $\SL(E)$-linearized locally free
sheaf $\Pr^k$ on $\SL(E)/P$. In this section we make the equivalence
explicit and give a geometric construction of $\Pr^k$ in terms of $\SL(E)$-linearized jet bundles on
$\flag(\dd,E)=\SL(E)/P$.

Let $\flag(\dd,E)=\SL(E)/P$ be the flag variety parametrizing flags of
type $\dd$ in an $N$-dimensional vector space $E$. Recall the Plucker
embedding
\[ i:\flag(\dd,E)\rightarrow \G=\G_1\times \cdots \times \G_k \]
defined by
\[ i([E_1\subseteq \cdots \subseteq E_k])=[E_1]\times \cdots \times
[E_k].\]
Let $q_i:\G\rightarrow \G_i$ be the projection morphism. Let
$\ll=(l_1,..,l_k)\in \mathbb{Z}^k$ and let
$\O(\ll)=q_1^*\O(l_1)\otimes \cdots \otimes q_k^*\O(l_k)$ be the
associated line bundle on $\G$. We get a linebundle
$\L(\ll)=i^*\O(\ll)$ on $\flag=\flag(\dd,E)$. Let $p,q:\flag \times
\flag\rightarrow \flag$ be the projection morphisms and let
$\I\subseteq \flag\times \flag$ be the ideal of the diagonal. 
\begin{definition} Let $\Pr^k_{\flag}(\L(\ll))=p_*(\O_{\flag\times
    \flag}/\I^{k+1}\otimes q^*\L(\ll))$ be the \emph{$k$'th order jet
    bundle} of $\L(\ll)$. 
\end{definition}

We first prove some general facts on jet bundles on arbitrary products
of schemes.

Let $A,B$ be arbitrary commutative $F$-algebras and let
$P^k_A=A\otimes_F A/I^{k+1}$ where $I\subseteq A\otimes_F A$ is the
ideal of the diagonal. If $X=\spec(A)$ it follows $\widetilde{P^k_A}=\Pr^k_X$.
There are natural maps of rings $p:A\rightarrow A\otimes B$ and
$q:B\rightarrow A\otimes B$.

\begin{lemma} \label{product}There is for every $k\geq 1$ a surjection of $A\otimes
  B$-modules
\[ P^k_A \otimes_F P^k_B\rightarrow P^k_{A\otimes B}.\]
\end{lemma}
\begin{proof} The natural map
\[ A\otimes A \otimes B \otimes B \rightarrow A\otimes B \otimes A
\otimes B \]
defined by
\[ \phi(a\otimes b \otimes x \otimes y)=a\otimes x \otimes b \otimes y
\]
indce a well defined map as claimed. It is surjective and the Lemma is proved.
\end{proof}

Let $A_1,..,A_s$ be commutative $F$-algebras.

\begin{corollary} \label{product2} There is for every $k\geq 1$ a surjective map of
  $A_1\otimes \cdots \otimes A_s$-modules
\[ P^k_{A_1}\otimes \cdots \otimes P^k_{A_s}\rightarrow
P^k_{A_1\otimes \cdots \otimes A_s}.\]
\end{corollary}
\begin{proof} The proof follows from Lemma \ref{product} and an induction.
\end{proof}

Let $E_i$ be an $A_i$-module for $i=1,..,s$.

\begin{corollary} \label{prod}There is for every $k\geq 1$ a surjection
\[ P^k_{A_1}(E_1)\otimes \cdots \otimes P^k_{A_s}(E_s)\rightarrow
P^k_{A_1\otimes \cdots \otimes A_s}(E_1\otimes \cdots \otimes E_s) \]
of $A_1\otimes \cdots \otimes A_s$-modules.
\end{corollary}
\begin{proof} This follows from directly Corollary \ref{product2}.
\end{proof}

Let $X_1,..,X_s$ be arbitrary schemes and let $\E_i$ be a quasi
coherent $\O_{X_i}$-module for $i=1,..,s$. Let $X=X_1\times \cdots
\times X_s$. Let $p_i:X\rightarrow X_i$ be the i'th projection and let 
$\E=p_1^*\E_1\otimes \cdots \otimes p_s^*\E_s$

\begin{corollary}\label{schemes} There is for every $k\geq 1$ a surjection
\[ p_1^*\Pr^k_{X_1}(\E_1)\otimes \cdots \otimes p_s^*\Pr^k_{X_s}(\E_s)\rightarrow
\Pr^k_{X}(\E) \]
of $\O_X$-modules.
\end{corollary}
\begin{proof} The Corollary is a global version of Corollary \ref{prod}.
\end{proof}

Let $q_i:\G\rightarrow \G_i$ be the projection morphism and let
\[ \O(\ll)=q_1^*\O(l_1)\otimes \cdots \otimes q_k^*\O(l_k) \]
be the line bundle on $\G$ defined above. 

\begin{proposition} Let $1\leq k \leq min\{l_i+1\}$. The $k$'th Taylor
  morphism
\[ T^k:\H^0(\G,\O(\ll))\rightarrow \Pr^k_{\G}(\O(\ll))(\overline{e})
\]
is surjective.
\end{proposition}
\begin{proof} There is by the Kunneth formula an isomorphism
\[ \H^0(\G,\O(\ll))\cong \H^0(\G_1,\O(l_1))\otimes \cdots
\otimes \H^0(\G_k,\O(l_k)) \]
of vector spaces.
The Taylor map $T^k_i$  is by \cite{verma} a surjective map
\[ T^k_i:\H^0(\G_i,\O(l_i))\rightarrow
\Pr^k_{\G_i}(\O(l_i))(\overline{e}) \]
for $i=1,..,k$.
We get a surjective map
\[\tilde{T}^k:\H^0(\G,\O(\ll))=\otimes^k_{i=1}
\H^0(\G_i,\O(l_i))\rightarrow^{\otimes_i T^k_i} \otimes_{i=1}^k
\Pr^k_{\G_i}(\O(l_i))(\overline{e}) \]
of vector spaces. By Corollary \ref{schemes} we get a surjective
morphism
\[ \otimes_{i=1}^k \Pr^k_{\G_i}(\O(l_i))(\overline{e})\rightarrow
\Pr^k_{\G}(\O(\ll)(\overline{e}) \]
of vector spaces.
This induce the surjection $T^k$ 
\[ T^k:\H^0(\G,\O(\ll))\rightarrow \Pr^k_{\G}(\O(\ll))(\overline{e})
\]
and the Proposition is proved.
\end{proof}

\begin{theorem} \label{surjective} Let $\L(\ll)\in \Pic^{\SL(E)}(\flag(\dd,E))$ be a line
  bundle with $l_i\geq 1$ for all $i$. Let $1\leq k \leq
  min\{l_i+1\}$. The Taylor map 
\[ T^k:\H^0(\flag(\dd,E),\L(\ll))\rightarrow
\Pr^k_{\flag}(\L(\ll))(\overline{e})\]
is a surjective map of vector spaces.
\end{theorem}
\begin{proof} Since $\L(\ll)=i^*\O(\ll)$ where $i:\flag \rightarrow \G$
  is the Plucker embedding, and the Taylor map is surjective on $\G$
  the Theorem follows from \cite{verma} Theorem 4.4.
\end{proof}

\begin{corollary} \label{jetexact}There is for $1\leq k \leq min\{l_i+1\}$ an exact
  sequence
\[ 0\rightarrow \H^0(\flag(\dd,E),\m^{k+1}\L(\ll))\rightarrow
\H^0(\flag(\dd,E),\L(\ll))\rightarrow
\Pr^k_{\flag}(\L(\ll))(\overline{e}) \rightarrow 0\]
of $P$-modules.
\end{corollary}
\begin{proof} Let $p,q:\flag\times \flag \rightarrow \flag$ be the
  projection morphisms and let $\I\subseteq \O_{\flag \times \flag}$
  be the ideal of the diagonal. Using higher direct images and the functor
  $p_*(-\otimes q^*\L(\ll))$ we get a long exact sequence of
  $\SL(E)$-linearized locally free sheaves
\[ 0\rightarrow p_*(\I^{k+1}\otimes q^*\L(\ll)) \rightarrow
p_*q^*\L(\ll)\rightarrow \Pr^k_{\flag}(\L(\ll))\rightarrow \]
\[ \R^1p_*(\I^{k+1}\otimes q^*\L(\ll))\rightarrow
\R^1p_*q^*\L(\ll)\rightarrow \cdots  .\]
Recall there is an equivalence of categories between the category of
$\SL(E)$ linearized vector bundles on $\flag(\dd,E)$ and the category
of rational $P$-modules.
We take the fiber a $\overline{e}\in \SL(E)/P$ to get an exact sequence
of $P$-modules
\[ 0\rightarrow \H^0(\flag, \m^{k+1}\L(\ll))\rightarrow
\H^0(\flag,\L(\ll)) \rightarrow^{T^k} \Pr^k_{\flag}(\L(\ll))(\overline{e})
\rightarrow \]
\[ \H^1(\flag,\m^{k+1}\L(\ll))\rightarrow \H^1(\flag,\L(\ll))\rightarrow
\cdots  \]
and since $\H^1(\flag,\L(\ll))=0$ and $T^k$ is surjective, the
Corollary follows.
\end{proof}

Dualize the exact sequence from Corollary \ref{jetexact} to get an
exact sequence of $\SL(E)$-modules
\[ 0\rightarrow \Pr^k_{\flag}(\L(\ll)(\overline{e})^*\rightarrow
\H^0(\flag,\L(\ll))^* \rightarrow^\phi
\H^0(\flag,\m^{k+1}\L(\ll))^*\rightarrow 0 .\]
The highest weight vector $v\in \H^0(\flag,\L(\ll))^*$ induce a
$P$-module
\[ \U^k(\lg)v\subseteq \H^0(\flag,\L(\ll))^* .\]

\begin{lemma} \label{inclusion}There is an inclusion of $P$-modules
\[ \U^k(\lg)v\subseteq \Pr^k_{\flag}(\L(\ll))(\overline{e})^* .\]
\end{lemma}
\begin{proof} One checks that $\phi(\U^k(\lg)v)=0$ and the Lemma follows.
\end{proof}

We can now prove the main theorem of the paper:
\begin{theorem} \label{maintheorem2} Let $1\leq k \leq min\{l_i+1\}$. There is an
  isomorphism of $P$-modules
\[ \U^k(\lg)v\cong \Pr^k_{\flag}(\L(\ll)(\overline{e})^* .\]
\end{theorem}
\begin{proof} There is by Lemma \ref{inclusion} an inclusion of
  $P$-modules
\[ \U^k(\lg)v\subseteq \Pr^k_{\flag}(\L(\ll))(\overline{e})^*.\]
By Corollary \ref{corrmain2} this inclusion is an isomorphism and the
Theorem is proved.
\end{proof}

We have for any finite dimensional irreducible $\SL(E)$-module $V$
with weight $\lambda$ constructed a linebundle $\L(\ll)$ on $\SL(E)/P$
where $P\subseteq \SL(E)$ is a parabolic subgroup with the following
property: There is an isomorphism $V\cong \H^0(\SL(E)/P,\L(\ll))^*$ of
$\SL(E)$-modules. Moreover the canonical filtration of $P$-modules
\begin{equation} \label{filt}
\U^1(\lg)v\subseteq \cdots \subseteq \U^k(\lg)v
\subseteq V
\end{equation}
equals the filtration
\begin{equation}\label{jetfilt}
 \Pr^1_{\flag}(\L(\ll))(\overline{e})^*\subseteq \cdots \subseteq
\Pr^k_{\flag}(\L(\ll))(\overline{e})^* \subseteq V=\H^0(\SL(E)/P,\L(\ll))^*
\end{equation}
given by the jet bundle $\Pr^i$. Here $1\leq k \leq min\{l_i+1\}$.
It follows Question \ref{main3} from the introduction is settled.

Assume $\underline{l}^i=(l^i_1,..,l^i_k)\in \Z^k$ with $l^i_j\geq 1$
for all $i,j$. Let $\E=\oplus_{i=1}^d \L(\underline{l}^i)$.
Let $v_i\in \H^0(G/P, \L(\underline{l}^i))$ be the unique highest
weight vector. Let $W\subseteq \H^0(G/P,\E)^*$ be the subspace
generated by $v_1,..,v_d$. Let 
\[ \U^l(\lg)W\subseteq \H^0(G/P,\E)^* \]
be the $P$-module generated by $W$ and $\U^l(\lg)$.

\begin{corollary} \label{reducible} There is an isomorphism
\[ \Pr^l(\E)(e)^*\cong \oplus_{i=1}^d \U(\lg)v_i \cong\U^l(\lg)W \]
of $P$-modules for all $1\leq l \leq min\{l^i_j+1\}$.
\end{corollary}
\begin{proof}
We get by Theorem \ref{maintheorem2} an isomorphism 
\[ \Pr^l(\E)(\oe)^*\cong \oplus_{i=1}^d
\Pr^l(\L(\underline{l}^i)(\oe)^* \cong  \oplus_{i=1}^d \U^l(\lg)v_i  \]
of $P$-modules,  and the claim of the Corollary follows.
\end{proof}

\begin{problem} Canonical filtrations for semi simple algebraic
  groups.\end{problem}

Let $G$ be any semi simple linear algebraic group over $F$ and let $V$
be any finite dimensional irreducible $G$-module with highest weight
vector $v\in V$. Let $L_v\subseteq V$ be the line spanned by $v$. Let
$P\subseteq G$ be the subgroup fixing the line $L_v$. It follows $P$
is a parabolic subgroup and the quotient $G/P$ is canonically a smooth
projective variety of finite type over $F$. Let $\lambda$ be the
weight of $v$ and let $\lg=Lie(G)$. Let $ann(v,\lambda)\subseteq
\U(g)$ be the left annihilator ideal of $v$ and let
$ann^k(v,\lambda)=ann^k(v,\lambda)\cap \U^k(\lg)$ be its canonical
filtration. We get an exact sequence
\[ 0\rightarrow ann(v,\lambda)\otimes L_v \rightarrow \U(\lg)\otimes
L_v \rightarrow V \rightarrow 0\]
of $G$-modules
and an exact sequence
\[0\rightarrow ann^k(v,\lambda)\otimes L_v \rightarrow
\U^k(\lg)\otimes L_v \rightarrow \U^k(\lg)L_v \rightarrow 0 \]
of $P$-modules.
\begin{definition} Let $\U^\bullet(\lg)L_v\subseteq V$ be the
  \emph{canonical filtration} of $V$.
\end{definition}

There is work in progress giving a geometric interpretation of the
canonical filtration $\U^\bullet(\lg)L_v$ in terms of $G$-linearized
$\O_{G/P}$-modules
$\Pr^k$ (see \cite{maa4}).

\begin{example} Morphisms of generalized Verma modules .\end{example}

Let $G$ be an arbitrary semi simple linear algebraic group and let
$P\subseteq G$ be a parabolic subgroup. Let $\lg=Lie(G)$ and
$\lp=Lie(P)$. Assume $U$ is a $G$-module and let $W,V\subseteq U$ be
sub $P$-modules. Assume
\[ f:\U(\lg)\otimes W\rightarrow \U(\lg)\otimes V \]
is a map of $G$-modules with
\begin{align} 
\label{a1}&f(\U^l(\lg)\otimes W )\subseteq \U^l(\lg)\otimes V \\
\label{a2}&f(ann(W)\otimes W)\subseteq ann(V)\otimes V.
\end{align}
The modules $\U(\lg)\otimes W$ and $\U(\lg)\otimes V$ are
\emph{generalized Verma modules}.
It follows $f$ induce a map
\[ f^l:\U^l(\lg)W\rightarrow \U^l(\lg)V \]
of $P$-modules. Here $\U^l(\lg)W$ and $\U^l(\lg)V$ are the
sub-$P$-modules generated by $\U^l(\lg)$, $W$ and $V$ as sub modules
of the $G$-module $U$. By the result of Corollary \ref{reducible} we
describe $\U^l(\lg)W$ in terms of $\Pr^l(\E)(\oe)^*$ for a locally
free $\O_{G/P}$-module $\E$ when $W$ is the
$P$-submodule generated by the highest weight vectors $v_i\in V_{\lambda_i}$
for $i=1,..,d$.

One seek to give a geometric construction of the
morphism $f^l$ in terms of $G$-linearized locally free
$\O_{G/P}$-modules:
We seek a morphism
\[ \phi:\E \rightarrow \F \]
of $G$-linearized $\O_{G/P}$-modules $\E,\F$ with
$\phi(\oe)=f^l$. This problem will be considered in later paper on
this subject (see \cite{matsumoto} for results on morphisms between
generalized Verma modules).

\textbf{Acknowledgements}: Thanks to Michel Brion, Alexei Roudakov 
and Claire Voisin for discussions and comments.

\end{document}